\documentclass[12pt,leqno]{compositio}
\usepackage{amsmath,amssymb,amscd,latexsym}

\usepackage[all]{xy}
\newcommand{\C}{{\mathbb{C}}}
\newcommand{\BP}{{\mathrm{BP}}}
\newcommand{\A}{{\mathrm{A}}}
\newcommand{\B}{{\mathrm{B}}}
\newcommand{\Cr}{{\mathrm{C}}}
\newcommand{\Zr}{{\mathrm{Z}}}
\newcommand{\MU}{{\mathrm{MU}}}
\newcommand{\Sr}{{\mathrm{S}}}
\newcommand{\D}{{\mathrm{D}}}
\newcommand{\Du}{{\underline{\underline{\mathrm{D}}}}}
\newcommand{\Er}{{\mathrm{E}}}
\newcommand{\Gm}{{\widehat{\mathbb{G}_m}}}
\newcommand{\F}{{\mathbb{F}}}
\newcommand{\Lr}{{\mathrm{L}}}
\newcommand{\Fr}{{\mathrm{F}}}
\newcommand{\M}{{\mathrm{M}}}

\newcommand{\oGamma}{\overline{\Gamma}}
\newcommand{\HH}{{\mathrm{H}}}

\newcommand{\Q}{{\mathbb{Q}}}
\newcommand{\Qr}{{\mathrm{Q}}}
\newcommand{\Z}{{\mathbb{Z}}}
\newcommand{\Zp}{{\mathbb{Z}_{(p)}}}

\newcommand{\id}{\mathrm{id}}
\newcommand{\imm}{\mathrm{im}\,}
\newcommand{\oT}{\overline{T}}
\renewcommand{\ker}{\mathrm{ker}}
\newcommand{\Spec}{\mathrm{Spec}\,}
\newcommand{\Ext}{\mathrm{Ext}}

\newcommand{\rk}{\mathrm{rk}\,}

\newcommand{\Fh}{{\mathcal F}}

\renewcommand{\oT}{\overline{T}}

\newtheorem{theorem}{Theorem}

\newtheorem{prop}[theorem]{Proposition}
\newtheorem{cor}[theorem]{Corollary}
\newtheorem{example}[theorem]{Example}
\newtheorem{remark}[theorem]{Remark}
\begin{document}

\title[Beta-elements and divided congruences]{Beta-elements and divided congruences}
\author{Jens Hornbostel and Niko Naumann}
\email{jens.hornbostel@mathematik.uni-regensburg.de\\
niko.naumann@mathematik.uni-regensburg.de}
\address{NWF I- Mathematik\\ Universit\"at Regensburg\\93040 Regensburg}

\begin{abstract}
The f-invariant is an injective homomorphism from the $2$-line of the Adams-Novikov spectral
sequence to a group which is closely related to divided congruences of elliptic modular forms.
We compute the f-invariant for two infinite families of $\beta$-elements and explain the
relation of the arithmetic of divided congruences with 
the Kervaire invariant one problem.
\end{abstract}

\maketitle

\section{Introduction}\label{intro}

One of the most successful tools for studying the stable
homotopy groups of spheres is the Adams-Novikov spectral sequence (ANSS)

\[ \Er_2^{s,t}=\Ext_{\MU_*\MU}^{s,t}(\MU_*,\MU_*)\Rightarrow 
\pi_{t-s}^s(\Sr^0).\]

The corresponding filtration on $\pi_*^s:=\pi_*^s(\Sr^0)$
defines a succession of invariants of framed bordism, each being defined whenever
all of its predecessors vanish, the first one of which is simply the degree

\[ d: \Fr^{0,*}/ \Fr^{1,*}=\pi_0^s\longrightarrow \Er_2^{0,0}=\Z \]

which is an isomorphism. The next invariant, defined for all $n>0$, is the 
e-invariant

\[ e: \pi_n^s=\Fr^ {1,n+1}\longrightarrow \Er_2^{1,n+1}\subseteq \Q/\Z,\]

c.f. \cite[Chapter 19]{Sw}. Though defined purely homotopy-theoretic here, the e-invariant
is well-known to encode subtle geometric information. For its relation to index theory
via the $\eta$-invariant see \cite[Theorem 4.14]{APS}. The e-invariant vanishes for all even $n=2k \ge 2$,
thus giving rise to the f-invariant

\[ f:\pi_{2k}^s=\Fr^{2,2k+2}\longrightarrow \Er_2^{2,2k+2}.\]

The understanding of this invariant is fragmentary at the moment. In particular, there is no
index-theoretic interpretation of it comparable to the one available for the e-invariant.\\
As a first step towards understanding the f-invariant, G. Laures \cite{L} showed how elliptic
homology can be used to consider the f-invariant 

\[ f: \pi_{2k}^s\longrightarrow \Er_2^{2,2k+2}\hookrightarrow \Du_{k+1}\otimes\Q/\Z \]

\medskip 

as taking values in a group which is closely related to divided congruences of 
modular forms. Note that this is similar to the role taken by complex $K$-theory
in the study of the e-invariant. Strictly speaking, at this point we had better switched from
$\MU$ to $\BP$. In fact, we will always work locally at a fixed prime $p$ in the following.\\
This surprising connection of stable homotopy theory with something 
as genuinely arithmetic
as divided congruences certainly motivates to ask for a thorough understanding of
how these are related by the f-invariant.\\
The main purpose of this paper is to make this relation explicit.\\
We also include a fairly self-contained review of G. Laures' above version 
of the f-invariant
to help the reader who might be interested in making his own computations. We now review
the individual sections in more detail.\\
In section \ref{computef}, we first remind the reader of the $\beta$-elements which generate
the $2$-line of the ANSS (with a little exception at the prime 2). We then construct, for 
suitable Hopf algebroids, a complex which is quasi-isomorphic to the cobar complex
and which will facilitate later computations. Finally, we show how to use elliptic
homology to obtain the f-invariant as above.\\
In section \ref{comp}, we recall some fundamental results of N. Katz on the arithmetic
of divided congruences and point out an interesting relation between $\BP$-theory and the
mod $p$ Igusa tower (Theorem \ref{igusa}). Next, we give some specific computations 
of modular forms and divided congruences for $\Gamma_1(3)$ which we will need to study the 
f-invariant of the Kervaire elements $\beta_{2^n,2^n}\in\Ext^{2,2^{n+2}}[\BP]$ at the prime $p=2$.\\ 
In section \ref{results}, we first compute the f-invariants of
the infinite family of $\beta$-elements $\beta_t$ for $t\ge 1$ not divisible by $p$ (Theorem \ref{betat}).
Then we explain how to approach the problem of
computing the Chern numbers determining the $\beta_{2^n,2^n}$ (see \cite[Corollary 4.2.5]{L2}
for the case of $\beta_1$ at the prime 3). We do this by explicit computations in $\BP$-theory
for dimension $2$ and $6$ (Theorem \ref{dastheorem}). The computations get very complicated in higher dimensions. In order to use divided congruences, we then compute the f-invariants
of the family $\beta_{s2^n,2^n}$ for $n\ge 0$ and $s\ge 1$ odd (Theorem \ref{fbeta}).
In Corollary \ref{detectkervaire} we determine
the divided congruence which detects the Kervaire element from stable homotopy.

\begin{acknowledgements} 
The idea that it might be possible 
to project to the element $\beta_{2^n,2^n}$ using divided 
congruences - and consequently 
rephrase the Kervaire invariant one problem - was communicated to us by G. 
Laures. We also thank the referee for a thorough report which 
led us to improve some points of the presentation.
\end{acknowledgements}

\section{The construction of the f-invariant}\label{computef}

We remind the reader of the construction of $\beta$-elements in section \ref{beta}.
In section \ref{ratcoeff}, we construct a complex which is quasi-isomorphic to the cobar complex
and in which we will compute representatives for some of the $\beta$-elements. This is used in
section
\ref{ell} where we explain how to express the f-invariant of elements 
in $\Ext^{2,2k}[\BP]$ in terms
of divided congruences.

\subsection{Beta-elements in stable homotopy}\label{beta}

We review some facts on Brown-Peterson homology $\BP$ at the prime $p$ 
and $\beta$-elements. See \cite{MRW} and \cite{R} for more details.
The Brown-Peterson spectrum $\BP$ has coefficient ring
$\BP_*=\Zp [v_1,v_2,...]$ with $v_i$ in dimension $2(p^i-1)$.
The universal $p$-typical formal group law is defined 
over this ring. The couple
$(\A,\Gamma):=(\BP_*,\BP_*\BP)$ becomes a Hopf algebroid
in a standard way and we have 
$\BP_*\BP=\BP_* \otimes \Zp [t_1,t_2,...]$ such that
the left unit $\eta_L$ of the Hopf algebroid
$(\A,\Gamma)$ is the standard inclusion. 
The right unit
$\eta_R$ is determined over $\Q$ by the formula in
\cite[Theorem A.2.1.27]{R}. Choosing the Hazewinkel generators
\cite[A.2.2.1]{R} for the $v_i$, a short computation
yields $\eta_R(v_1)=v_1+ p t_1$ and $\eta_R(v_2)
=v_2 + v_1t_1^p - v_1^pt_1$ mod $p$.\\
We have the chromatic resolution of $\BP_*$ as a left $\BP_*\BP$-comodule
$$ \BP_* \to M^0 \to M^1 \to ... $$
which gives rise to the chromatic spectral sequence
$$\Ext^{s,*}_{\BP_*\BP}(\BP_*,M^t) \Rightarrow \Ext^{s+t,*}_{\BP_*\BP}(\BP_*,\BP_*).$$
This allows the construction of elements in $\Ext^{*,*}_{\BP_*\BP}(\BP_*,\BP_*)$,
the so called Greek-letter elements. Strictly speaking, these elements
arise from comodule sequences $0 \to N^n \to M^n \to N^{n+1} \to 0$,
but for our computations we will need the related comodule sequences
(\ref{1}) and (\ref{2}) below. See \cite[Lemma 3.7 and Remark 3.8]{MRW}
for the relationship between them. 
We abbreviate $\HH^n(\cdot):=\Ext^{n,*}_{\BP_*\BP}(\BP_*,\cdot)$ in the following.\\
To construct the $\beta$-elements \cite[p. 476/477]{MRW}, choose integers $t,s,r\geq 1$ 
such that $(p^r,v_1^s,x_n^{t'})\subseteq \BP_*$ is an invariant ideal where $t=p^nt'$, $(p,t')=1$
and $x_n$ is a homogeneous polynomial in $v_1, v_2$ and $v_3$
considered as an element of $v_2^{-1}\BP_*/(p^r,v_1^s)$ (see 
\cite[p. 202]{R} or \cite[p. 476]{MRW}), for example $x_0=v_2$.\\
Consider the two short exact sequences of $\BP_*\BP$-comodules

\begin{equation}\label{1}
 0 \to \BP_* \stackrel{p^r}{\to} \BP_* \to \BP_*/(p^r) \to 0
\end{equation}

\begin{equation}\label{2}
0 \to \Sigma^{2s(p-1)} \BP_*/(p^r) \stackrel{v_1^s}{\to} \BP_*/(p^r) \to \BP_*/(p^r,v_1^s) \to 0.
\end{equation}
 
Using the induced boundary maps 
\[ \delta: \HH^0(\BP_*/(p^r,v_1^s))\longrightarrow \HH^1(\BP_*/(p^r))\mbox{ and }\]

\[ \delta': \HH^1(\BP_*/(p^r))\longrightarrow \HH^2(\BP_*) \]

we define $\beta_{t,s,r}:=\delta' \delta(x_n^{t'})$.\\
It is known \cite[Theorem 2.6]{MRW}, \cite{Sh} for which indices $(t,s,r)$
the elements $ \beta_{t,s,r}$ are non-zero in $\HH^2(BP_*)$. In this case the order
of $\beta_{t,s,r}$ is $p^r$. By construction, we have $\beta_{t,s,r}\in
\HH^{2,2t(p^2-1)-2s(p-1)}(\BP_*)$.\\

\begin{example} 
For $p=2$ and $n \geq 1$, the element $\beta_{2^n,2^n}:=\beta_{2^n,2^n,1} \in 
\Ext^{2,2^{n+2}}_{\BP_*\BP}(\BP_*,\BP_*)$
is called the {\it Kervaire element}. It is mapped via the Thom reduction
\cite[Theorem 5.4.6]{R} to the element $h_{n+1}^2 \in 
\Ext^{2,2^{n+2}}_{\HH \Z /2_* \HH \Z /2}(\HH \Z /2_*,\HH \Z /2_*)$.
The latter element survives to a non-zero element of $\pi_{2^{n+2}-2}(\Sr^0)$ in the 
Adams spectral sequence if and only if (as W. Browder \cite[Theorem 7.1]{Br}
has shown)
there exists a framed manifold of dimension $2^{n+2}-2$
with non-vanishing Kervaire invariant. Whether or not this is the case is unknown for
$n \geq 5$, for $n\leq 4$ see \cite{BJM}, \cite{KM}.
\end{example}

\subsection{The rationalised cobar complex}\label{ratcoeff}

The standard way of displaying elements in $\Ext^n:=\Ext^n_{\Gamma}(\A,\A)$ of a flat Hopf algebroid $(\A,\Gamma)$ is
by means of the cobar complex. In this section we shall give another 
description of this Ext group as a subquotient of $(\A\otimes \Q)^{\otimes n}$
needed to compute f-invariants (Proposition \ref{rat}). The results of this section are a more algebraic version of
\cite[Section 3.1]{L}.\\
Let $(\A,\Gamma)$ be a Hopf algebroid with structure maps $\eta_L,\eta_R,\epsilon$ and $\Delta$.
This determines a cosimplicial abelian group $\Gamma^{\cdot}$ as follows:
Set $\Gamma^n:=\Gamma^{\otimes_{\A} n}$ with cofaces
 $\partial^i:\Gamma^n\longrightarrow\Gamma^{n+1}$;
 $\partial^0(\gamma_1\otimes\ldots\otimes\gamma_n)
:=1\otimes\gamma_1\otimes\ldots\otimes\gamma_n$;
 $\partial^i(\gamma_1\otimes\ldots\otimes\gamma_n)
:=\gamma_1\otimes\ldots\otimes\Delta(\gamma_i)
\otimes\ldots\otimes\gamma_n$ ($1\leq i \leq n$)
 and $\partial^{n+1}(\gamma_1\otimes\ldots
\otimes\gamma_n):=\gamma_1\otimes\ldots\otimes\gamma_n
\otimes 1$
 for $n\geq 1$ and $\partial^0:=\eta_R$,
 $\partial^1:=\eta_L$ for $n=0$ and
codegeneracies $\sigma^i:\Gamma^{n+1}\longrightarrow
\Gamma^n$,
 $\sigma^i(\gamma_0\otimes\ldots\otimes\gamma_n):=
\gamma_0\otimes\ldots\otimes\epsilon(\gamma_i)\otimes
\ldots\otimes\gamma_n$.
We also denote by $\Gamma^{\cdot}$ the
associated cochain complex.\\
Following \cite[Definition A.1.2.11]{R}, we define
the reduced cobar complex (usually denoted as $\Cr_{\Gamma}^{\cdot}(\A,\A)$) as being the 
subcomplex $\overline{\Gamma}^{\cdot}\subseteq\Gamma^{\cdot}$ with $\overline{\Gamma}^n:=\overline{\Gamma}^{\otimes_{\A} n}$ for $n\geq 1$ where $\overline{\Gamma}:=\ker(\epsilon)$ and $\overline{\Gamma}^0:=\A$. This is
a subcomplex because $\Delta(\gamma)-\gamma\otimes 1- 1 \otimes \gamma\in\overline{\Gamma}\otimes_{\A}\overline{\Gamma}$ for any $\gamma\in\overline{\Gamma}$. 

We now assume that $(\A,\Gamma)$ is a flat Hopf algebroid such that \\
i) $\A$ (and hence $\Gamma$) is torsion free.\\
ii) The map $\A_{\Q}^{\otimes 2}:=(\A\otimes \Q)^{\otimes 2}\stackrel{\phi}{\longrightarrow} \Gamma_{\Q}:=\Gamma\otimes \Q$, $a\otimes b\mapsto a.\eta_R(b)$ is an isomorphism ($\eta_L$ is suppressed from notation).\\
iii) The $\Q$-algebra $\A_{\Q}$ is augmented by some $\tau:\A_{\Q}\longrightarrow\Q$.\\

\begin{remark}
The above assumptions i)-iii) are fulfilled by the flat Hopf algebroid $(\BP_*,\BP_*\BP)$: 
The proof of ii) follows from the fact that over any $\Q$-algebra any two
($p$-typical) formal group laws are isomorphic via a unique strict isomorphism. 
If $\BP_*\longrightarrow \Er_*$ is a non-zero Landweber exact algebra \cite{HS}, 
then $(\Er_*,\Er_*\Er)$ also fulfils all the above assumptions: $\Er_*$ is a $\Z_{(p)}$-algebra on which $p$ is not a zero divisor, hence $\Er_*$ is torsion free.
Conditions ii), iii) and the flatness of of $(\Er_*,\Er_*\Er)$ are inherited from
$(\BP_*,\BP_*\BP)$ because $\Er_*\Er\simeq \Er_*\otimes_{\BP_*} \BP_*\BP \otimes_{\BP_*} \Er_*$ 
, c.f. \cite[Proposition 10]{N} for the flatness.
\end{remark}

We define the cosimplicial abelian group $\D^{\cdot}$ by $\D^n:=\A_{\Q}^{\otimes (n+1)}$ ($n\geq 0$) with cofaces $\partial^i:\D^n\longrightarrow \D^{n+1}$ to be given by
$\partial^i(a_0\otimes\ldots\otimes a_n):=a_0\otimes\ldots\otimes 1\otimes\ldots\otimes a_n$
with the $1$ in the $(i+1)^{st}$ position ($i=0,\ldots,n+1$) and codegeneracies $\sigma^i:\D^n\longrightarrow \D^{n-1}$ ($i=0,\ldots,,n-1)$ defined by $\sigma^i(a_0\otimes\ldots\otimes a_n):=a_0\otimes\ldots\otimes a_ia_{i+1}\otimes\ldots\otimes a_n$. By ii) above
we have for any $n\geq 0$ an isomorphism

\begin{equation}\label{isoDG}
 \phi^n:\D^n=\A_{\Q}^{\otimes(n+1)}\simeq(\A_{\Q} \otimes_{\Q} \A_{\Q})^{\otimes_{\A_{\Q}}n}\stackrel{\phi^{\otimes n}}{\longrightarrow}\Gamma_{\Q}^{\otimes_{\A_{\Q}}n}=\Gamma^n \otimes \Q
\end{equation}

which maps $a_0\otimes\ldots\otimes a_n\mapsto a_0\otimes\ldots\otimes a_{n-2}\otimes a_{n-1}.\eta_R(a_n)$ and one
checks that $\phi^{\cdot}$ is an isomorphism of cosimplicial groups and hence of cochain complexes.\\
We have $\HH^{\cdot}(\D^{\cdot})=\HH^0(\D^{\cdot})=\Q$, a contracting homotopy being given by

\begin{equation}\label{contract}
 H^n: \D^n=\A_{\Q}^{\otimes(n+1)}\longrightarrow \D^{n-1}=\A_{\Q}^{\otimes n}\; ;
a_0\otimes\ldots\otimes a_n\mapsto \tau(a_0)a_1\otimes\ldots\otimes a_n. 
\end{equation}

Define a subcomplex $\Sigma^{\cdot} \subseteq \D^{\cdot}$ by

\[ \Sigma^n:=\sum_{i=1}^n\partial^i(\D^{n-1})=\sum_{i=1}^n \A_{\Q}^{\otimes i} \otimes \Q \otimes \A_{\Q}^{\otimes (n-i)}\subseteq \D^n=\A_{\Q}^{\otimes (n+1)},\]

for $n\geq 1$ and $\Sigma^0:=0$. One checks that the composition

\[ \iota^n: \overline{\Gamma}^n\hookrightarrow \Gamma^n \hookrightarrow \Gamma^n\otimes\Q\stackrel{(\phi^n)^{-1}}{\longrightarrow} \D^n\longrightarrow \D^n/\Sigma ^n \]
is injective for all $n\geq 0$. This is obvious for $n=0$, and for $n\geq 1$ it follows from
$(\cap_{i=0}^{n-1}\ker(\sigma^i))\cap(\sum_{i=1}^n \imm(\partial^i))=0$, which in turn is 
an easy consequence of the cosimplicial identities: One
shows by descending induction on $1\leq j\leq n$ that $(\cap_{i=0}^{n-1}\ker(\sigma^i))
\cap(\sum_{i=j}^n \imm(\partial^i))=0$. Observe that $\D^n/\Sigma ^n$
is isomorphic to the group labelled $E_* \otimes G_*^n$ in \cite[p. 404]{L} for $\A=\BP$.
\\
We define a cochain complex $\Qr^{\cdot}$ 
by the exactness of

\begin{equation}\label{defQ}
0\longrightarrow \overline{\Gamma}^{\cdot} \stackrel{\iota^{\cdot}}{\longrightarrow} \D^{\cdot}/\Sigma^{\cdot}\stackrel{\pi^{\cdot}}{\longrightarrow}\Qr^{\cdot}\longrightarrow 0,
\end{equation}

hence $\Qr^n\simeq \A_{\Q}^{\otimes (n+1)}/\Sigma^n+\imm(\overline{\Gamma}^{\otimes_\A n})$.\\
From the definitions of the differential of $\D^{\cdot}$ and $\Sigma^n$ one obtains
($\B\Qr^n\subseteq \Qr^n$ denoting the boundaries)

\[ \Qr^n/\B\Qr^n\simeq \A_{\Q}^{\otimes (n+1)}/\tilde{\Sigma}^n+\imm(\overline{\Gamma}^{\otimes_{\A} n}), \]

where $\tilde{\Sigma}^n:= \Sigma^n + \Q\otimes \A_{\Q}^{\otimes n} =
\sum_{i=1}^{n+1} \A_{\Q}^{\otimes (i-1)}\otimes\Q \otimes \A_{\Q}^{\otimes (n+1-i)}$. 

The alternative description of $\Ext$ we are aiming for is the following.

\begin{prop}\label{rat}
For any $n\geq 1$, the connecting homomorphism $\delta$ of (\ref{defQ}) 
is an isomorphism 

\[ \HH^n(\Qr^{\cdot})\stackrel{\delta}{\longrightarrow} \HH^{n+1}(\overline{\Gamma}^{\cdot})= \Ext^{n+1}. \]
\end{prop}

\begin{proof} One readily sees that the contracting homotopy (\ref{contract})
of $\D^{\cdot}$ respects the subcomplex $\Sigma^\cdot$ in positive dimensions,
hence the middle term of (\ref{defQ}) is acyclic in these dimensions.
\end{proof}

To explicitly compute $\delta$, it is useful to note that the differential of
$\D^{\cdot}/\Sigma^{\cdot}$ has the simple form

\[ \D^n/\Sigma^n\stackrel{d}{\longrightarrow} \D^{n+1}/\Sigma^{n+1}\;,\;
[a_0\otimes\ldots\otimes a_n]\mapsto [1\otimes a_0\otimes\ldots\otimes a_n], \]

as is immediate from the definitions. To compute $\delta^{-1}$, we consider the zig-zag

\[ \xymatrix{ \Ext^{n+1} \ar@{.>}[rrrrrd]^-{\delta^{-1}} & \Zr(\overline{\Gamma}^{n+1}) \ar[r]^-{\iota^{n+1}} \ar@{>>}[l] & \D^{n+1}/\Sigma^{n+1} \ar[r]^-{H^{n+1}} & \D^n/\Sigma^n \ar@{>>}[r]^-{\pi^n} & \Qr^n \ar@{>>}[r] & \Qr^n/\B\Qr^n \\
 & & & & & \HH^n(\Qr^{\cdot}). \ar@{_{(}->}[u]  } \]

One checks that the dotted arrow exists and is the inverse of $\delta$.\\
Now let $p$ be a prime, $r\geq 1$ an integer, and let $\delta':\Ext^n_{\Gamma}(\A,\A/p^r)\longrightarrow \Ext^{n+1}$
be the connecting homomorphism associated to the short exact sequence of 
$\Gamma$-comodules $0\longrightarrow \A\stackrel{\cdot p^r}{\longrightarrow} \A \longrightarrow \A/p^r \longrightarrow 0$. Consider the diagram

\begin{equation}\label{spare}
\xymatrix{ \delta':\Ext^n_{\Gamma}(\A,\A/p^r) & \Zr(\overline{\Gamma}^n\otimes_{\A} \A/p^r) \ar@{>>}[l] & \{ z\in \overline{\Gamma}^n | dz\in p^r\overline{\Gamma}^{n+1} \}  \ar@{>>}[l] \ar[r]^-{\alpha} \ar[d]_{\nu} \ar@{.>}[rd] & \Ext^{n+1} \\
 & & \Qr^n/\B\Qr^n &  \HH^n(\Qr^{\cdot}) \ar@{_{(}->}[l] \ar[u]_{\simeq}^{\delta}. }
\end{equation}

Here, $\alpha(z):=[y]$ for any $y\in \overline{\Gamma}^{n+1}$ satisfying $dz=p^ry$ and
$\nu(z):=\pi^n(p^{-r}\iota^n(z))$ mod $\B\Qr^n$. The upper horizontal line
is $\delta'$ by definition and one checks that $\nu$ factors through
the dotted arrow and makes the diagram commutative.
Hence, when displaying an element of $\imm(\delta')$ 
in $\HH^{\cdot}(\Qr^{\cdot})$ rather than in the usual cobar complex, one does not have to compute 
the cobar differential implicit in $\alpha$ but only the (easier) map $\nu$.\\
Finally, assume that everything in sight is graded where the grading on
$\D^n$,$\Gamma^n$, etc. is by total degree. For a fixed $k\in\Z$, consider
the commutative diagram 

\[ \xymatrix{ \Qr^n/\B\Qr^n \cong \A_{\Q}^{\otimes(n+1)}/\tilde{\Sigma}^n+\imm(\overline{\Gamma}^{\otimes_{\A} n}) & \HH^n(\Qr^{\cdot}) \ar@{_{(}->}[l] \ar[r]^{\delta}_{\simeq} & \Ext^{n+1} \\
\A_{\Q}^{\otimes (n+1),k}/\tilde{\Sigma}^{n,(k)}+\imm(\overline{\Gamma}^{\otimes_{\A} n})^k \ar@{^{(}->}[u]_j & \HH^{n,k}(\Qr^{\cdot}) \ar@{^{(}->}[u] \ar[r]^{\simeq} \ar@{_{(}->}[l] & \Ext^{n+1,k}. \ar@{_{(}->}[u]  } \]

Here $\tilde{\Sigma}^{n,(k)}:=\tilde{\Sigma}^n\cap \A_{\Q}^{\otimes (n+1),k}$.
One checks that $j$ is well defined and injective, and $\HH^{n,k}(\Qr^{\cdot})$ is defined to
be the pull back of $\HH^n(\Qr^{\cdot})$ along $j$. The commutative diagram
induces an isomorphism $\HH^{n,k}(\Qr^{\cdot})\stackrel{\simeq}{\longrightarrow}\Ext^{n+1,k}$ as indicated.\\
For example, for $(\A,\Gamma)=(\BP_*,\BP_*\BP)$ and $n=1$ we obtain an inclusion

\begin{equation}\label{inclusion}
 \Ext^{2,k}\subseteq\frac{(\BP_{\Q}\otimes \BP_{\Q})^{(k)}}{\BP_k\BP+(\BP_{\Q}\otimes\Q+\Q\otimes \BP_{\Q})^{(k)}} 
\end{equation}

which is important for us since the f-invariant is defined in terms of the 
group on the right hand side.\\
To effectively compute representatives of $\beta$-elements in the complex $\Qr^{\cdot}$ one proceeds as follows. Let $t,s,r \geq 1$ be integers as in section \ref{beta} and $\delta$, $\delta '$ the coboundary maps introduced there. Fix $k\in\Z$ and $x\in \HH^{0,k}(\A/(p^r,v_1^s))$, that is

\[ x\in \Cr_{\Gamma}^{0,k}(\A/(p^r,v_1^s))=(\A/(p^r,v_1^s))^k \]

is an invariant element ($\Cr_{\Gamma}$ indicates the reduced cobar complex).
As $\delta$ is the connecting homomorphism determined by 
the short exact sequence of complexes obtained by applying
$\Cr_{\Gamma}$ to (\ref{2}), we compute $\delta(x)$ as
follows: Lift $x$ to $y\in (\A/(p^r))^k$ and compute the cobar differential

\[ d=\eta_R-\eta_L: \Cr_{\Gamma}^{0,k}(\A/p^r)=(\A/p^r)^k\longrightarrow
(\A/p^r\otimes_{\A}\oGamma)^k=(\oGamma/p^r)^k=\Cr_{\Gamma}^{1,k}(\A/p^r) \]

obtaining $d(y)\in(\oGamma/p^r)^k=\oGamma^k/p^r$.
Note that this computation requires knowledge of $\eta_R(y)$ mod $p^r$.\\
Now $d(y)$ will be divisible by $v_1^s$, hence

\[ d(y)=v_1^sz\mbox{ in } (\oGamma/p^r)^k \]

with $z\in (\oGamma/p^r)^{k-2s(p-1)}=\Cr_{\Gamma}^{1,k-2s(p-1)}(\A/p^r)$
representing $\delta(x)$. Lift $z$ to some $w\in\oGamma^{k-2s(p-1)}$.
To proceed, we use diagram (\ref{spare}): $w$ lies in the second to last group
of the top row, hence we compute $\nu(w)\in \HH^2(Q^\cdot)$ which is our 
representative for $\delta'(\delta(x))\in \Ext^2$. This requires to compute
$\phi^{-1}$. Observe for example that $\phi^{-1}(t_1)=
\frac{1\otimes v_1 - v_1 \otimes 1}{p}$.\\

\subsection{Elliptic homology theories and divided congruences: the f-invariant}\label{ell}

We refer the reader to \cite{L} or \cite{H} for the notion of elliptic homology 
with respect to the congruence subgroup $\Gamma_1(N)$. In this 
section, $\Er$ denotes the spectrum associated to the homology theory with
coefficient ring $\Er_*=\M_*(\Z_{(p)},\Gamma_1(N))[\Delta^{-1}]$, see section \ref{forms} for the 
notation. Finally, $\alpha: \BP \longrightarrow \Er$ denotes the orientation.\\
By the naturality of the constructions in section \ref{ratcoeff} we have a commutative diagram for any $k\ge 0$

\[ \xymatrix{ \Ext^{2,k}[\BP] \ar@{^{(}->}[r]^{\alpha} \ar@{^{(}->}[d]_{(\ref{inclusion})} & \Ext^{2,k}[\Er] \ar@{^{(}->}[d]^{\iota}_{(\ref{inclusion})} \\
\frac{(\BP_{\Q}\otimes \BP_{\Q})^{(k)}}{\BP_k\BP+(\BP_{\Q}\otimes\Q+\Q\otimes \BP_{\Q})^{(k)}} \ar[r]^{\alpha\otimes\alpha} & \frac{(\Er_{\Q}\otimes \Er_{\Q})^{(k)}}{\Er_k\Er+(\Er_{\Q}\otimes\Q+\Q\otimes \Er_{\Q})^{(k)}} \,.} \]

The injectivity of $\alpha$ holds for any Landweber exact theory $\Er$ of height at least two, \cite[Proof of 4.3.2]{L2}.
To proceed, however, we will use a more subtle property of $\Er$, namely
the topological $q$-expansion principle. We put 

\[ \Du_k:=\{ f=\sum_{i=0}^kf_i\in\bigoplus_{i=0}^k\Er_{\Q,2i} | \mbox{there are } g_0\in\Q,g_k\in\Er_{\Q,2k}\mbox{ such that } (f+g_0+g_k)(q)\in\Zp^{\Gamma}[[q]]\} \]

where $f(q)$ denotes the $q$-expansion of $f$ at the cusp infinity,
and for $\Gamma=\Gamma_1(N)$ we set
$\Z^{\Gamma}:=\Z[\frac{1}{N},\zeta_N]$ if $N>1$ and
$\Z^{SL_2(\Z)}:=\Z[\frac{1}{6}]$ as in \cite{L}.
We then define

\begin{equation}\label{iota2}
 \iota^2: \frac{(\Er_{\Q}\otimes \Er_{\Q})^{(2k)}}{\Er_{2k}\Er+(\Er_{\Q}\otimes\Q+\Q\otimes \Er_{\Q})^{(2k)}}\longrightarrow\Du_k\otimes\Q/\Z\, , \, \sum_{i+j=k}f_i\otimes g_j\mapsto \sum_{i+j=k}-q^0(f_i)g_j\, ,
\end{equation}

where $q^0(f)$ is the constant term of the $q$-expansion of $f$ at the cusp infinity. The composition $\iota^2\circ\iota$ is injective \cite[Proposition 3.9]{L} and hence 
so is the f-invariant
\[ f:\Ext^{2,2k}[\BP]\hookrightarrow\Du_k\otimes\Q/\Z\,.\]

We remark that our grading of $\Er_*$ is the topological one, i.e. elements of dimension $2k$
correspond to modular forms of weight $k$.

G. Laures describes $\Ext^2$ using the 
 canonical Adams resolution \cite[Definition 2.2.10]{R}
instead of the cobar resolution. 

\begin{prop}
For $k > 0$ even, the above definition of the $f$-invariant $\Ext^{2,k} \longrightarrow \Du_{k/2}\otimes \Q/\Z$
coincides with the one given in \cite{L}.
\end{prop}
\begin{proof}
We have to show that the map 
$\Ext^{2,k} \hookrightarrow \frac{(\BP_{\Q}\otimes \BP_{\Q})^{(k)}}
{\BP_k\BP+(\BP_{\Q}\otimes\Q+\Q\otimes \BP_{\Q})^{(k)}}$ we constructed
in section \ref{ratcoeff},(\ref{inclusion}) coincides
with the one of \cite{L}.
We know \cite[Lemma A.1.2.9 (b)]{R} that, up to chain homotopy, there is a unique map from the unreduced cobar resolution
$\Gamma \otimes_{\A} \Gamma^{\otimes_{\A} \cdot}$ to the
canonical Adams resolution $\BP_*(\BP \land \Sigma \bar{\BP}^{\land \cdot})$
where $\bar{\BP} \to S^0 \stackrel{\eta}{\to} \BP \stackrel{d}{\to} \Sigma \bar{\BP}$
is an exact triangle in the stable homotopy category, see also \cite[Lemma 3.7]{Bru}.
Now one checks that $\Gamma ^{\otimes_{\A} (n+1)} \cong \pi_*(\BP^{\land n+2})
\stackrel{\pi_*(id_{\BP^{\land 2}} \land d^{\land n})}{\to}
\pi_*(\BP^{\land 2}\land(\Sigma\bar{\BP})^{\land n})$ is a map of chain complexes where
the isomorphism follows from \cite[Lemma 2.2.7]{R} and induction.
So the claim reduces to the fact that the triangle
$$ \xymatrix{ \Zr^{2,k}(\Gamma^{\otimes_{\A} \cdot}) \ar[rr]^{\pi_*(\id_{\BP^{\land 2}} \land d)} 
\ar[dr] & & 
\Zr^{2,k}(\pi_*(\BP^{\land 2} \land (\Sigma \bar{\BP})^{\land \cdot})) \ar[dl] \\
 & \frac{(\BP_{\Q}\otimes \BP_{\Q})^{(k)}}
{\BP_k\BP+(\BP_{\Q}\otimes\Q+\Q\otimes \BP_{\Q})^{(k)}} & } $$
commutes. This follows using that our maps 
$\tau, H^2$ and $\pi_*(d \land d)$
correspond to $r, \rho$ and the isomorphism
$\D^1/\tilde{\Sigma^1} \cong G^2= \pi_*((\Sigma \bar{\BP})^{\land 2}) \otimes \Q$ 
in \cite[section 3.1]{L} (where we define $\tau$ by mapping all
$v_i$ to zero).
\end{proof}

Note the degree shift and the factor 2 for the f-invariant
$\pi_{2k}^s\longrightarrow \Du_{k+1}\otimes\Q/\Z$ (both are missing
in \cite[p. 411]{L}).

\section{Arithmetic computations}\label{comp}

In section \ref{div}, we review results of N. Katz on divided congruences and establish a relation
between $\BP$-theory and the mod $p$ Igusa tower (Theorem \ref{igusa}). In section \ref{forms}
we give explicit computations for elliptic homology of level 3 and the corresponding divided congruences.

\subsection{Divided congruences}\label{div}

We review parts of \cite{K}. Some technical remarks are in order: In {\em loc. cit.} N. Katz
works with level-$N$ structures of fixed determinant for $N\ge 3$,
denoted by $\Gamma^{can}(N)$ in \cite{KaM}. To confirm
with general policy in algebraic topology we wish to consider $\Gamma_1(N)$-structures instead,
which are representable only for $N\ge 5$. Note that the results of \cite{K} carry over to this moduli problem since the irreducibility of the Igusa tower is independent of level structures, c.f. \cite[Theorem 4.3 and Theorem 4.3 bis]{K2}. Furthermore, we will use 
these results for $p\ge5$ and $N=1$. These cases can be handled by using
auxiliary rigid level structures and taking invariants under suitable finite groups as in 
\cite{K}.\\
Fix a level $N\ge 5$, a prime $p$ not dividing $N$ and a primitive $N$-th root of 
unity $\zeta\in\overline{\F_p}$. We put $k:=\F_p(\zeta)$, $W:=W(k)$ (Witt vectors) and also
denote by $\zeta\in W$ the Teichm\"uller lift of $\zeta$. Finally, $K$ denotes the
field of fractions of $W$ and for any $\Z_{(p)}$-algebra $R$ we denote by $\M_k(R,\Gamma_1(N))$
the $R$-module of holomorphic modular forms for $\Gamma_1(N)$ of weight $k$ and defined over $R$ , see e.g. \cite{K2} or \cite{L}.
The $\Gamma_1(N)$ is omitted from the notation if it is clear from the context.
We fix a lift $E_{p-1}\in\M_{p-1}(W)$ of the Hasse invariant.
The existence of such a lift puts further restrictions on $p$ and $N$ which are satisfied in our applications. That is, we consider either the case
$p \geq 5$ (see \cite[p. 98]{K2}) or provide an explicit lift
(as after Remark 8). \\
We define the ring of divided congruences $\D$ by

\[ \M_*(W)\subseteq \D:=\{ f\in\M_*(K) | f(q)\in W[[q]] \} \subseteq \M_*(K), \]

where $f(q)$ is the $q$-expansion of $f$ at the cusp infinity. For $n\ge 0$ we also define

\[ \D_n:=\D\cap \left( \bigoplus_{i=0}^n \M_i(K)\right)\subseteq\D\mbox{ and} \]
\[ \Du_n:=\D_n+K+\M_n(K)\subseteq  \bigoplus_{i=0}^n \M_i(K) \]

which is consistent with the definition of the previous section.
The group $\underline{\underline{D}}_k$
considered in $\cite{L}$ differs from the $\Du_k$ above because the ring
of holomorphic modular forms has been localised in \cite{L}.
This difference is not serious because the $f$-invariant
factors through holomorphic modular forms
\cite[Proposition 3.13]{L}.

The ring $\D$ carries a uniformly continuous $\Z_p^*$-action (the diamond operators) defined by

\[ [\alpha](\sum_i f_i):=\sum_i \alpha^if_i,\]

where $\alpha\in\Z_p^*$ and $f_i\in\M_i(K)$. We put $\Gamma_0:=\Z_p^*$ and
$\Gamma_n:=1+p^n\Z_p$ for $n\ge 1$ and $V_{1,n}:=(\D/p\D)^{\Gamma_n}$ for $n\ge 0$. Then

\[ V_{1,0}\subseteq V_{1,1}\subseteq\ldots\subseteq \D/p\D \]

is an ind-\'etale $\Z_p^*$-Galois extension, the mod $p$ Igusa tower. So, $V_{1,0}\subseteq V_{1,1}$
is a $(\Z/p)^*$-Galois extension and for all $n\ge 2$ $V_{1,n-1}\subseteq V_{1,n}$ is an \'etale
$\Z/p$-extension and hence an Artin-Schreier extension. An immediate computation with 
diamond operators shows that the composition $\M_*(W)\to\D\to\D/p\D$ factors through $V_{1,1}\subseteq \D/p\D$.
It is a result of P. Swinnerton-Dyer that this induces an isomorphism

\begin{equation}\label{PSD}
 \M_*(W)/(p,E_{p-1}-1)\stackrel{\simeq}{\longrightarrow} V_{1,1},
\end{equation}

see \cite[Corollary 2.2.8]{K}.\\
By Artin-Schreier theory, given $n\ge 2$ and $x\in\D/p\D$ satisfying $[\alpha](x)=x$ for 
all $\alpha\in\Gamma_n$ and $[1+p^{n-1}](x)=x+1$ one has $V_{1,n}=V_{1,n-1}[x]$ and the minimal 
polynomial of $x$ over $V_{1,n-1}$ is $T^p-T-a$ for some $a\in V_{1,n-1}$.\\
At this point we can establish a first relation between $\BP$-theory and divided congruences.
Consider the ring extensions 

\[ \M_*(W)\subseteq\D\subseteq W[[q]]. \]

We have a formal group $\Fh$ over $\M_*(W)$ induced by the universal elliptic curve. The base change
of $\Fh$ to $W[[q]]$ is the formal completion of a Tate elliptic curve and is thus isomorphic
to $\Gm$. Implicit in \cite{K} is the fact that $\D$ is the {\em minimal} extension of $\M_*(W)$
over which $\Fh$ becomes isomorphic to $\Gm$, i.e. $\D$ is obtained from $\M_*(W)$ by adjoining
the coefficients of an isomorphism $\Fh\simeq\Gm$ defined over $W[[q]]$. This is what underlies 
N. Katz' construction \cite[Section 5]{K} of a sequence of elements $d_n\in\D$ which modulo $p$
constitute a sequence of Artin-Schreier generators for the mod $p$ Igusa tower.\\
Since the elements $t_n\in\BP_*\BP$ are the coefficients of the universal isomorphism of a 
$p$-typical formal group law, one may expect a relation between the $t_n$ and the $d_n$. To formulate
this, denote by $\alpha:\BP_*\longrightarrow \M_*(W)$ the classifying map of $\Fh$ and consider the composition

\[ \rho: \BP_*\BP\subseteq\BP_*\BP\otimes\Q\stackrel{(\ref{isoDG})}{\simeq}\BP_{\Q}^{\otimes 2}\stackrel{\alpha\otimes\alpha}{\longrightarrow}\M_*(K)^{\otimes 2}\stackrel{-q^0\otimes\id}{\longrightarrow}\M_*(K). \]

Note that the map $\iota_2$ in (\ref{iota2}) composed with the orientation
$\alpha$ is a quotient of $\rho$. \\
The topological $q$-expansion principle guarantees that $T_n:=\rho(t_n)\in\D$ for $n\ge 1$ and
we can thus define $\oT_n:=(T_n$ mod $p\D)\in\D/p\D$.

\begin{theorem}\label{igusa}

For any $n\ge 1$ we have $[1+p^k](\oT_n)=\oT_n$ for $k>n$ and $[1+p^n](\oT_n)=\oT_n+1$.\\
Hence $\oT_n$ is an Artin-Schreier generator for the extension $V_{1,n}\subseteq V_{1,n+1}$.
\end{theorem}

\begin{proof} Let $\omega=(\sum_{n\ge 1} a_nt^{n-1})dt$ be the expansion along infinity
of a normalised (i.e. $a_1=1$) invariant differential on the universal elliptic curve.
Then $a_n\in\M_{n-1}(W)$ and the logarithm of the $p$-typification is $\sum_{n\ge 0}\frac{a_{p^n}}{p^n}t^{p^n}\in\M_*(K)[[t]]$, i.e. the classifying map $\alpha:\BP_*\longrightarrow \M_*(W)$,
when tensored with $\Q$, sends $l_n\in\BP_{\Q,2(p^n-1)}$ to $\frac{a_{p^n}}{p^n}\in\M_{p^n-1}(K)$,
see \cite[Theorem A.2.1.27]{R} for the definition of the $l_n$ (=$\lambda_n$ in the notation of
{\em loc. cit.}).\\
Defining $d_0:=1$ and $d_n$ ($n\ge 1$) recursively by

\begin{equation}\label{defdn}
\sum_{i=0}^n\frac{d_{n-i}^{p^i}}{p^i}=\frac{a_{p^n}}{p^n},
\end{equation}

N. Katz shows in \cite[Corollary 5.7]{K} that the $\overline{d}_n:=(d_n$ mod $p)\in\D/p\D$ behave under
the diamond operators as claimed for the $\oT_n$. In $\BP_*\BP\otimes\Q$ we have
$\eta_R(l_n)=\sum_{i=0}^nl_it_{n-i}^{p^i}$ and we apply $\rho$ to this relation to obtain

\begin{equation*}
\frac{a_{p^n}}{p^n}=\sum_{i=0}^n\frac{q^0(a_{p^i})}{p^i}T_{n-i}^{p^i}
\end{equation*}

which, using (\ref{defdn}), implies

\begin{equation}\label{dT}
\sum_{i=0}^n\frac{d_{n-i}^{p^i}}{p^i}=\sum_{i=0}^n\frac{q^0(a_{p^i})}{p^i}T_{n-i}^{p^i}.
\end{equation}

We now proceed by induction on $n\ge 1$. For $n=1$ we have $d_1+1/p=q^0(a_1)T_1+q^0(a_p)/p$.
Also, $q^0(a_1)=1$ since $a_1=1$ and $q^0(a_p)\in 1+pW$ because $a_p$ reduces mod $p$ to the Hasse
invariant which has $q$-expansion equal to $1$. We obtain $\oT_1=\overline{d}_1+\alpha$ for some
$\alpha\in k$. As $\alpha$ is invariant under all diamond operators, our claim for $n=1$ is
obvious.\\
Assume that $n\ge 2$. From (\ref{dT}) and $a_1=1$ we obtain 

\[ T_n=\sum_{i=0}^n\frac{d_{n-i}^{p^i}}{p^i}-\sum_{i=1}^n\frac{q^0(a_{p^i})}{p^i}T_{n-i}^{p^i}. \]
For $k>n$ we know that the terms involving $d_i$ are invariant mod $p$ under $[1+p^k]$ whereas
the remaining terms are likewise invariant by the induction hypothesis.\\
Finally, we have 

\[ [1+p^n]T_n=[1+p^n]d^n+[1+p^n](\sum_{i=1}^n\frac{d_{n-i}^{p^i}}{p^i})-[1+p^n](\sum_{i=1}^n\frac{q^0(a_{p^i})}{p^i}T_{n-i}^{p^i}). \]

Here we have $[1+p^n]d^n\equiv d_n+1\,(pD)$ and the remaining terms are invariant. Thus, indeed,
$[1+p^n](\oT_n)=\oT_n+1$.
\end{proof}

\subsection{Modular forms}\label{forms}

For a prime $p\ge 5$, the following is well known \cite[Appendix]{L}: 

\[ \M_*(\Z_{(p)},\Gamma_1(1))=\Z_{(p)}[E_4,E_6],\]

where $E_4$ and $E_6$ are the Eisenstein series of level one of the indicated weight.
For the discriminant $\Delta$, the ring of meromorphic modular forms is given by $\Z_{(p)}[E_4,E_6,\Delta^{-1}]$ and the usual orientation 

\[ \BP_*\longrightarrow \Z_{(p)}[E_4,E_6,\Delta^{-1}] \]

is Landweber exact of height $2$ and factors through $\Z_{(p)}[E_4,E_6]$.
A similar result holds for $p \geq 3$ and 
$\M_*(\Z_{(p)},\Gamma_1(2))=\Zp[\delta, \epsilon]$. \\
The purpose of this section is to give analogous results for $\Gamma_1(3)$ and $p=2$, c.f. \cite{St}
for related results.\\
Consider the elliptic curve

\[ E: y^2+a_1xy+a_3y=x^3 \]

defined over $R:=\Z[1/3][a_1,a_3,\Delta^{-1}]$ where $\Delta=a_3^3(a_1^3-27a_3)$ is the discriminant
of the given Weierstrass equation. Note that, unlike in level one, $\Delta$ is not irreducible
as a polynomial in $a_1$ and $a_3$ and we put $f:=a_3$, $g:=a_1^3-27a_3$, hence $\Delta=f^3g$.\\
The section $P:=(0,0)\in E(R)$ is of exact order 3 in every geometric fibre as follows from
\cite[III,2.3]{Si1} and $\omega:=dx/(2y+a_1x+a_3)$ is an invariant differential on $E$.

The following may be compared with \cite[Lemma 11]{St}:

\begin{prop}\label{universal}
The above tuple $(E/R,\omega,P)$ is the universal example of an elliptic curve over a $\Z[1/3]$-scheme together with a point of order $3$ and a non-zero invariant differential.
\end{prop}

\begin{proof} We have to show that whenever $T$ is a $\Z[1/3]$-scheme and $E'/T$ is an elliptic 
curve with non-zero invariant differential $\omega'$ and $P'\in E'(T)$ of exact order $3$, there is
a unique map $\phi:T\longrightarrow\Spec(R)$ such that $\phi^*(E,P,\omega)=(E',P',\omega')$.
We show the uniqueness of $\phi$ first. This amounts to seeing that the only change of coordinates

\[ x=u^2x'+r\, , \, y=u^3y'+u^2sx'+t \]

with $r,s,t\in R$ and $u\in R^*$ (see \cite[III Table 1.2]{Si1}) preserving $(E,P,\omega)$
is the identity, i.e. $r=s=t=0$ and $u=1$.\\
From $x'(P)=y'(P)=0$ we obtain $r=t=0$. Next, $a_4=a_4'$ implies $-sa_3=0$, hence $s=0$
because $\Delta=f^3g$ and thus also $f=a_3$ is a unit in $R$. Finally,
$\omega'=u\omega$ forces $u=1$.\\
Given the uniqueness of $\phi$ in general, its existence is a local problem on $T$ and we can
assume that $T=\Spec(S)$ is affine and $E'/T$ is given by a Weierstrass equation with coefficients
$a'_i\in S$. Moving $P'$ to $(0,0)$ gives $a_6'=0$. We claim that $a'_3\in S^*$: 
This can be checked on geometric fibres where it follows from \cite[III,2.3]{Si1} and the fact that
$(0,0)$ has order $3$ (if $a'_3$ vanished on some geometric fibre the point $(0,0)$ would have
order $2$ in that fibre). Using this, one finds a transformation such that $(dy)_{P'}=0$ in
$\Omega_{E'/T,P'}$, hence $a'_2=a'_4=0$. We thus have some $\psi:T\longrightarrow \Spec(R)$ 
such that $\psi^*(E,P)=(E',P')$ and $\psi^*(\omega)=u\omega'$ for some $u\in S^*$.
Adjusting $\psi$ using $u$, i.e. multiplying the $a'_i$ by $u^{-i}$, we obtain the desired $\phi$.
\end{proof}

We conclude that the ring of meromorphic modular forms is given as

\[ \M_*^{mer}(\Z_{(2)},\Gamma_1(3))=\Z_{(2)}[a_1,a_3,\Delta^{-1}] \]

with $a_i$ of weight $i$, and likewise for any other prime different from $3$ in place of $2$.\\
As usual, $t=-x/y$ is a local parameter at infinity for $E/R$ which is normalised for $\omega$
and hence determines a $2$-typical formal group law over $\M_*^{mer}(\Z_{(2)},\Gamma_1(3))$.
Using \cite[p. 113]{Si1} one checks that the corresponding classifying map 

\[ \alpha: \BP_*\longrightarrow \M_*^{mer}(\Z_{(2)},\Gamma_1(3)) \]

satisfies $\alpha(v_1)=a_1$ and $\alpha(v_2)=a_3$ for the Hazewinkel generators $v_i$.
Thus $\alpha$ makes $\M_*^{mer}(\Z_{(2)},\Gamma_1(3))$ a Landweber exact $\BP$ algebra of height $2$.\\
Using the orders of $f$ and $g$ at the two cusps $0$ and $\infty$ of $X_1(3)$, one can check 
that the ring of holomorphic modular forms is given by 

\begin{equation}\label{holo}
 \M_*(\Z_{(2)},\Gamma_1(3))=\Z_{(2)}[a_1,a_3]. 
\end{equation}

We stick to the notations of section \ref{div} for $p=2$ and $N=3$. For example, $\zeta$
denotes a primitive cube-root of unity and $W=W(\F_2(\zeta))=W(\F_4)=\Z_2[\zeta]$ is the unique unramified quadratic extension of $\Z_2$.\\
To study divided congruences we will need to know the $q$-expansions of $a_1$ and $a_3$. 
Given a Dirichlet character $\chi$, we consider it as a function on $\Z$ as
usual and define for $k\ge 0$ and $n\ge 1$

\[ \sigma_k^{\chi}(n):=\sum_{1\leq d|n}\chi(d)d^k. \]

In the following, $\chi$ will always denote the unique non-trivial character mod $3$

\[ \chi: (\Z/3\Z)^*\longrightarrow\C^*.\]

\begin{prop}\label{qexp}
The $q$-expansions of $a_1$ and $a_3$ at the cusp infinity are given as follows.
\[ a_1(q)=(1+2\zeta)(1+6\sum_{n\ge 1}\sigma_0^{\chi}(n)q^n)\mbox{ and }\]
\[ a_3(q)=(1+2\zeta)(-\frac{1}{9}+\sum_{n\ge 1}\sigma_2^{\chi}(n)q^n)\mbox{ in }W[[q]].\]
\end{prop}

\begin{proof} From (\ref{holo}) we know that $\rk\M_1(\Z_{(2)})=1$ and $\rk\M_3(\Z_{(2)})=2$.
In \cite[section 2.1.1]{Kp} we find the following modular forms:

\begin{equation}\label{katz}
6G_{1,\chi}(q)=1+6\sum_{n\ge 1}\sigma_0^{\chi}(n)q^n\in\M_1(\Z_{(2)})\mbox{ and}
\end{equation}

\[ G_{3,\chi}(q)=-\frac{1}{9}+\sum_{n\ge 1}\sigma_2^{\chi}(n)q^n\in\M_3(\Z_{(2)}),\]

where we have evaluated $\Lr(0,\chi)=1/3$ and $\Lr(-2,\chi)=-2/9$ using \cite[Theorem VII.2.9]{Ne} and \cite[formula following Proposition 4.1 and Exercise 4.2(b)]{Wa}.\\
It is easy to see that $G_{3,\chi}(0)=0$, i.e. $G_{3,\chi}$ vanishes at the cusp $0$: Note that $G_{3,\chi}$ is an Eisenstein series (c.f. Remark \ref{Hecke} below) and use equations (4) and (2) of \cite{He}.
Below, we will explain how to compute the following values of $a_1$ and $a_3$ at the cusps zero
and infinity.
\begin{equation}\label{x}
a_1(\infty)=1+2\zeta
\end{equation}
\begin{equation}\label{y}
a_3(\infty)=-\frac{1}{9}(1+2\zeta)
\end{equation}
\begin{equation}\label{z}
a_3(0)=0.
\end{equation}

Using these values and the dimensions of the spaces of
modular forms of weight $1$ and $3$, we conclude that $a_1=6(1+2\zeta)G_{1,\chi}$ and 
$a_3=(1+2\zeta)G_{3,\chi}$, hence that $a_1$ and $a_3$ have the desired $q$-expansions by (\ref{katz}). We are using the fact that the map $\M_3(\Z_{(2)})\otimes\C=\M_3(\C)\longrightarrow\C^2\, ,\, f\mapsto (f(\infty),f(0))$ is an isomorphism, as follows from the theory of 
Eisenstein series.\\
To establish $(\ref{x})$ and $(\ref{y})$ one has to evaluate $a_1$ and $a_3$ at the 
tuple $(T(q),\omega_{can},P)$ consisting of the Tate curve $T(q)/\Z((q))$, its canonical invariant
differential $\omega_{can}$ and a specific section $P\in T(q)(\Z[\zeta]((q)))[3]$. To do so,
one may use J. Tate's uniformisation \cite[p. 426]{Si2} to write $T(q)/(\Z[[q]]/(q^3))$
in Weierstrass form, the point $P$ having coordinates $(X(q,\zeta),Y(q,\zeta))$. One then uses Weierstrass
transformations to bring $(T(q),\omega_{can},P)/\Z[[q]]/(q^3)$ to the standard form of
Proposition \ref{universal}. The coefficients $a_1$ and $a_3$ of the Weierstrass equation thus obtained are by definition $a_1(\infty)$ and $a_3(\infty)$. The computation for (\ref{z}) is
similar, the point $P$ has to be replaced by $Q=(X(q,q^{1/3}),Y(q,q^{1/3}))$.
\end{proof}

\begin{remark}\label{Hecke}
In E. Hecke's notation \cite{He}, we have $a_1=\frac{9i}{\pi}G_1(\tau,0,1,3)$
and $a_3=\frac{27i}{4\pi^3}G_3(\tau,0,1,3)$.
\end{remark}

Note that $a_1(q)\equiv 1$ mod $2$, hence $a_1\in\M_1(\Z_{(2)},\Gamma_1(3))$ is a lift of the Hasse invariant for $p=2$.\\
Observe further that the results
of section \ref{div} carry over to the moduli
problem $\Gamma^{can}(N)$ initially studied
by Katz in \cite{K}, including the case $p=2$ and $N=3$.
We have an inclusion of rings of modular forms
$\M_*(W,\Gamma_1(3)) \subseteq \M_*(W,\Gamma^{can}(3))$
which induces a morphism of rings of divided 
congruences $D^{\Gamma_1(3)} \to D^{\Gamma^{can}(3)}$. 
The subgroups of $D^{\Gamma^{can}(3)}/2D^{\Gamma^{can}(3)}$ in the Igusa tower 
are denoted by $V_{1,n}^{can}$.
From section \ref{div}, we also know that 
\[ V_{1,0}^{can}=V_{1,1}^{can}\stackrel{(\ref{PSD})}{\simeq}
\M_*(W,\Gamma^{can}(3))/(2,a_1-1) \supseteq
\M_*(W,\Gamma_1(3))/(2,a_1-1)=k[a_3] \]
($k:=W/2W=\F_4$) and that, for $T:=\frac{q^0(a_1)-a_1}{2}\in\D$, $\oT:=(T$ mod $2)
\in\D^{\Gamma^{can}(3)}/2\D^{\Gamma^{can}(3)}$
is an Artin-Schreier generator for $V_{1,1}^{can}\subseteq V_{1,2}^{can}$. In particular $\oT^2+\oT\in V_{1,1}^{can} \supseteq k[a_3]$ and for later use we will need the following more precise result.

\begin{prop}\label{relation} $\oT^2+\oT=1+a_3$.
\end{prop}
\begin{proof} Recall that the $q$-expansion map $V_{1,1}^{can}\subseteq D^{\Gamma^{can}(3)}/2 D^{\Gamma^{can}(3)} \hookrightarrow k[[q]]$ is injective \cite[(1.4.6) for $m=1$]{K}. In $k[[q]]$ we have $T=\sum_{n\ge 1}\sigma_0^{\chi}(n)q^n$ and 
$T^2=\sum_{n\ge 1}\sigma_0^{\chi}(n)q^{2n}$, hence

\[ T^2+T=\sum_{n\ge 1}(\sigma_0^{\chi}(n/2)+\sigma_0^{\chi}(n))q^n, \]

where we understand that $\sigma_0^{\chi}(n/2)=0$ for $n$ odd.
To complete the proof, one needs to check that for all $n\ge 1$ one has

\[ \sigma_0^{\chi}(n/2)+\sigma_0^{\chi}(n)\equiv \sigma_2^{\chi}(n)\mbox{ mod }2, \]

and we leave this exercise in elementary number theory to the reader.
\end{proof}

\section{f-invariants and Kervaire invariant one}\label{results}

In this section, we compute the f-invariants of two infinite families of $\beta$-elements including the Kervaire elements $\beta_{2^n,2^n}$ and explain the relation of our results with
the Kervaire invariant one problem.

\subsection{ $f(\beta_t)$ for $t$ not divisible by $p$}\label{beta1}

Fix a prime $p$ and the level $N$ as $N=1$ for $p\ge 5$ and $N=5$ for $p=2,3$. 
We keep the notations of section \ref{div} for this choice of $p$ and $N$. Given an integer $t\ge 1$ not divisible by $p$, recall that 
$\beta_t=\delta' \delta (v_2^t)\in\Ext^{2,2t(p^2-1)-2(p-1)}
[\BP]$ has its f-invariant 
in $\Du_n\otimes\Q/\Z$, $n:=t(p^2-1)-(p-1)$.\\
When trying to express $f(\beta_t)$ in terms of divided congruences, we encounter what is
in fact the major obstacle at the moment for using the arithmetic of divided congruences in
homotopy theory:\\
The group $\Du_n\otimes\Q/\Z$ is not {\em directly} related to $\D$.
Instead, we have 
$\Du_n=\D+K+\M_n(K)$ by definition and there is a canonical surjection

\[ \pi:\D_n\otimes\Q/\Z\simeq\left(\bigoplus_{i=0}^n\M_i(K)\right)/\D_n \longrightarrow \Du_n\otimes\Q/\Z
\simeq  \left(\bigoplus_{i=0}^n\M_i(K)\right)/\Du_n\]

which is split because its kernel is divisible, hence $W$-injective.
In particular, $\pi$ remains surjective when restricted to $p$-torsion

\begin{equation}\label{quot}
\pi: \D_n\otimes\Q/\Z[p] \longrightarrow \Du_n\otimes\Q/\Z[p],
\end{equation}

note that $f(\beta_t)\in \Du_n\otimes\Q/\Z[p]$. The group $\D_n\otimes\Q/\Z[p]$ is related
to the ring of divided congruences as follows:

\begin{equation}\label{mappy}
\psi: \D_n\otimes\Q/\Z[p]\stackrel{\simeq}{\longrightarrow}\D_n/p\D_n\hookrightarrow\D/p\D,
\end{equation}

where the first arrow is multiplication by $p$ and the injectivity of the last map
is immediate. What we will do is to compute some element in $\D/p\D$ in the image of 
$\psi$ which under $\pi$ projects to $f(\beta_t)$. At the low risk of confusion we 
will continue to label such an element, which is in general not unique, as $f(\beta_t)$.
Recall that we have fixed an elliptic orientation $\alpha:\BP_*\longrightarrow\M_*(W)$
and denote $T:=\frac{\alpha(v_1)-q^0(\alpha(v_1))}{p}\in\D/p\D$. We also put $b:=((q^0(\alpha(v_2))$ mod $p)\in k$.

\begin{theorem}\label{betat}
For an integer $t\ge 1$ not divisible by the fixed prime $p$, we have
\[ f(\beta_t)=b^t-(T^p-T+b)^t\in V_{1,0}\subseteq\D/p\D.\]
\end{theorem}

\begin{proof} Note first that from section \ref{div} we know that $T$ is an Artin-Schreier generator
for $V_{1,1}\subseteq V_{1,2}$, hence $T^p-T\in V_{1,1}$. A short computation with diamond operators, which we leave to the reader, shows that in fact $T^p-T\in V_{1,0}$, hence
also $b^t-(T^p-T+b)^t\in V_{1,0}$.\\
We introduce $a:=q^0(\phi(v_1))$ and compute as explained at the end of section \ref{ratcoeff} using the notations 
introduced there. From $\eta_R v_2\equiv v_2+v_1t_1^p-v_1^pt_1$ mod $p$ we obtain
\[ \eta_R v_2^t\equiv v_2^t+\sum_{i=1}^t { t \choose i} v_2^{t-i}v_1^it_1^i(t_1^{p-1}-v_1^{p-1})^i\mbox{ mod } p ,\]
hence
\[ w= \sum_{i=1}^t { t \choose i} v_2^{t-i}v_1^{i-1}t_1^i(t_1^{p-1}-v_1^{p-1})^i\]
and 
\[ \nu(w)=\frac{1}{p}\sum_{i=1}^t { t \choose i}(v_2^{t-i}v_1^{i-1}\otimes 1)\left( \frac{1\otimes v_1 - v_1 \otimes 1}{p} \right) ^i\cdot \]
\[ \left( \left( \frac{1\otimes v_1 - v_1 \otimes 1}{p}\right) ^{p-1}- v_1^{p-1}\otimes 1 \right) ^i\in\frac{(\BP_{\Q}\otimes \BP_{\Q})^{(2n)}}{\BP_{2n}\BP+(\BP_{\Q}\otimes\Q+\Q\otimes \BP_{\Q})^{(2n)}}.\]

As in section \ref{ell} we apply $\iota^2\circ(\alpha\otimes\alpha)$ to this expression to 
obtain, denoting $\alpha(v_1)\in\M_{p-1}(W)$ as $v_1$ for simplicity,

\[ -\frac{1}{p} \sum_{i=1}^t {t \choose i} b^{t-i}a^{i-1}\left( \frac{v_1-a}{p}\right) ^i \left( \left( \frac{v_1-a}{p}\right) ^{p-1} - a^{p-1} \right) ^i =\]
\[ \frac{-1}{pa}\left( \left( a\left( \frac{v_1-a}{p}\right) \left( \left( \frac{v_1-a}{p} \right) ^{p-1} - a^{p-1}\right) +b \right)^t -b^t\right )\in\bigoplus_{i=0}^n\M_i(K).\]
This is a representative for $f(\beta_t)$ in $\D_n\otimes\Q/\Z[p]$ to which we have to apply the map $\psi$ from (\ref{mappy})
to obtain an element in $\D/p\D$. For this, note that $a\equiv 1$ mod $p$ because $v_1$
reduces to the Hasse invariant mod $p$. This allows us to put $a=1$ in the above expression
(but not to replace $\frac{v_1-a}{p}$ by $\frac{v_1-1}{p}$; this would require the congruence
$a\equiv 1$ mod $p^2$, which does not hold in general). We then obtain indeed

\[ -((T(T^{p-1}-1)+b)^t-b^t)=b^t-(T^p-T+b)^t.\]

\end{proof}

\begin{remark}\label{laures} Assume that $p\ge 5$ in the situation of Theorem \ref{betat}. In general, the elliptic orientation will not map $v_1$ to the Eisenstein series $E_{p-1}$ of weight $p-1$ and level one. But $\alpha(v_1)$ and $E_{p-1}$ can only differ by a modular form divisible by $p$ and we may thus change the orientation to force $\alpha(v_1)=E_{p-1}$. Assuming this, we see that $f(\beta_{1,1,1})= 
\frac{E_{p-1}-1}{p^2} - \frac{1}{p}
(\frac{E_{p-1}-1}{p})^p$, as first computed by G. Laures \cite[p. 414]{L}
(where the second summand is missing).
\end{remark}

\begin{remark} The injectivity of the f-invariant together with the known structure
of $\Ext^2[\BP]$ provides some non-trivial information about the arithmetic of divided congruences as follows. 
Fix some $x\in\Ext^{2,k}[\BP]$ of order $p^r$. Then $f(x)\in \Du_{k}\otimes\Q/\Z$ will be of order
$p^r$, hence a representative of $f(x)$ in $\D_k\otimes\Q/\Z$ will be of order $p^s$ for some
$s\ge r$. Thus the f-invariant relates the order of a $\beta$-element 
to the (non-)existence of a certain divided congruence. \\
Let us assume that $r=1$ as is the case for all $\beta_{t,s,r}$
considered in this article. Then our results show that our representatives 
in $\D\otimes\Q/\Z$ have order $p$ and the non-trivial
additional information on divided congruences is then that they do not lie in the kernel of $\D\otimes\Q/\Z\longrightarrow\Du\otimes\Q/\Z$.\\
To give an example, assume that we are in the situation of Remark \ref{laures}. 
The arithmetic of divided congruences shows that $F:=\frac{E_{p-1}-1}{p^2} - \frac{1}{p}
(\frac{E_{p-1}-1}{p})^p\in\D_{p(p-1)}\otimes\Q/\Z$ is of order $p$, i.e. $F\in\bigoplus_{i=0}^{p(p-1)}\M_i(K)$ has a $q$-expansion with denominator exactly $p$. The additional information is then that for any $\alpha\in K$ and $f\in\M_{p(p-1)}(K)$ the $q$-expansion of $F+\alpha+f$ will still have exact denominator $p$.
\end{remark}

\begin{example}
Fix $p=5$ and set $g_2:=\frac{1}{12}E_4$ and $g_3:=\frac{-1}{216}E_6$
as in \cite{K}.
The comparison of the logarithm of the universal $p$-typical formal group law 
\cite{R} and the corresponding coefficients of the logarithm of
the elliptic curve $(E, \omega)$ \cite[(5.0.3)]{K} (p-typification does not 
change these coefficients)
shows that the orientation $a$ maps $v_1$ to $a_p$ and
$v_2$ to $\frac{a_{p^2}-a_p^{p+1}}{p}$, the $a_i$ denoting
the normalised (multiplied with $- 1/2$) $a_i$
of \cite[p. 351]{K}. 
One deduces that $v_1$ maps to $-8g_2$ and $v_2$ maps to 
$\frac{a_{25}-a_5^6}{5}$. 
A computation with Maple 
shows that $a_{25}=129761280 g_2^3g_3^2 + 32440320 g_3^4 + 3784704 g_2^6$ 
(and also that the correct value for the unnormalised
$a_{11}$ is $-2520g_2g_3$ and not $-512g_2g_3$).
It follows that $q^0(v_1)=-\frac{2}{3}$ and $q^0(v_2)=-\frac{4900}{3^{10}}$,
so Theorem \ref{betat} may be rephrased in terms of the Eisenstein series
$g_2$ and $g_3$.
\end{example}

\subsection{ Projecting to the Kervaire element}\label{kervaire}

In this section, we compute $f(\beta_{s2^n,2^n})$ for $n\ge 0$ and $s\ge 1$ odd at the prime $p=2$.
Using this, we are able to determine a single coefficient in the f-invariant of a $(U,fr)^2$-
manifold of dimension $2^n$ the non-vanishing of which is necessary and sufficient
for the corner of $X$ to be a Kervaire manifold, that is having Kervaire
invariant one. See \cite{L2} for the notion
of cobordism of manifolds with corners.\\
We begin by recalling the well-known relation of the Kervaire invariant to certain $\beta$-elements,
due to W. Browder \cite{Br}. Fix some $n\ge 3$. We have a homomorphism

\[ K:\pi^s_{2^n-2}\longrightarrow\Z/2 \]

which sends the class of a stably framed manifold to its Kervaire invariant. Consider on
the other hand the composition 

\[ K':\pi^s_{2^n-2}\longrightarrow\pi^s_{2^n-2}[2]\longrightarrow \Er_{\infty}^{2,2^n}[\HH\Z/2] \hookrightarrow\Er_2^{2,2^n}[\HH\Z/2]=\Z/2\cdot h_{n-1}^2. \]

Here, the first map is the projection to the $2$-primary part, the second is the projection onto
$F^2/F^3$ in the (classical) Adams spectral sequence at $p=2$, the third is an edge homomorphism and
the final equality is due to J. Adams, \cite[3.4.1, c)]{R}.

\begin{prop}\label{Kprime} $K=K'$.
\end{prop}

\begin{proof} For any $y\in\pi^s_{2^n-2}[2]$ we have $K(y)=1$ if and only if $y$ has Adams filtration $2$. This is implicit in \cite{Br}, c.f. \cite[p. 144]{bjm}.
\end{proof}

We can easily obtain a similar homotopy theoretic description of $K$ using $\BP$
instead of $\HH\Z/2$. 

\begin{prop}\label{shimo} 
Let $n\geq 2$. Then $\Ext^{2,2^{n}}$ is a direct sum of cyclic groups of order
$2$.  It is generated by the element
$\alpha_1 \cdot {\alpha_{2^{n-1}-1}}$ and
the elements $\beta_{s2^i,2^i}$ with $s$ odd and $i \geq 0$ 
such that $(3s-1)2^{i+1}=2^n$ and the case $(s,i)=(1,0)$ has to be omitted.
\end{prop}
\begin{proof}
This follows from \cite[Corollary 5.4.5]{R}.
Observe that the $\overline{\alpha}_t$ in loc. cit 
equals $\alpha_t$ as $t$ is odd, see \cite[Theorem 1.5]{Sh}
or \cite[Theorem 5.2.6]{R}. 
\end{proof}

\begin{remark} The Lemma shows that the number of generators of $\Ext^{2,2^{n}}$ is 
$[n/2]+1$ for $n\ge 3$. The low dimensional cases are as follows.
\[ Ext^{2,4}\;:\; \alpha_1^2 \]
\[ Ext^{2,8}\;:\; \alpha_1 {\alpha}_3,\beta_{2,2} \]
\[ Ext^{2,16}\;:\; \alpha_1 {\alpha}_7,\beta_{4,4},\beta_{3,1} \]
\[ Ext^{2,32}\;:\; \alpha_1 {\alpha}_{15}, \beta_{8,8}, \beta_{6,2} \]
\[ Ext^{2,64}\;:\; \alpha_1 {\alpha}_{31}, \beta_{16,16}, \beta_{12,4}, 
\beta_{11,1}\]
\[ Ext^{2,128}\;:\; \alpha_1 {\alpha}_{63}, \beta_{32,32}, \beta_{24,8},
\beta_{22,2} \]
\[ Ext^{2,256}\;:\; \alpha_1 {\alpha}_{127}, \beta_{64,64}, \beta_{48,16},
\beta_{44,4}, \beta_{43,1}.\]
\end{remark}

Now we consider the composition

\[ K'':\pi^s_{2^n-2}\longrightarrow \pi^s_{2^n-2}[2]\longrightarrow \Er_{\infty}^{2,2^n}[\BP]\hookrightarrow\Er_2^{2,2^n}[\BP]\longrightarrow\Z/2\cdot\beta_{2^{n-2},2^{n-2}} \]

which is defined in analogy with $K'$, the final map being the projection to the $\Z/2$-summand generated by $\beta_{2^{n-2},2^{n-2}}$.

\begin{prop}\label{kervaireBP} $K''=K$.
\end{prop}

\begin{proof} We have the Thom reduction $\Phi:\Ext^*[\BP]\longrightarrow
\Ext^*[\HH\Z/2]$ which satisfies
$\Phi(\beta_{2^{n-2},2^{n-2}})=h_{n-1}^2$ and is zero on all other generators of $\Ext^{2,2^n}[\BP]$, see \cite[5.4.6, a)]{R} and Proposition \ref{shimo}. The result then 
follows from Proposition \ref{Kprime}.
\end{proof}

Let $X$ be a $(U,fr)^2$-manifold of dimension $2^n$. From the above, we see that the
corner of $X$ is a Kervaire manifold if and only if the f-invariant of $X$ contains
$\beta_{2^n,2^n}$ as a summand. Thus, one certainly wants a more geometric description of
the f-invariant (or just its projection to $\beta_{2^{n-2},2^{n-2}}$). In principle, it is possible to obtain such a description in terms of Chern numbers of $X$, simply because they determine the
$(U,fr)^2$-bordism class of $X$ \cite{L2}, but the necessary computations become quite 
complicated already
in low dimensions. The next theorem gives a sample of such a computation.
At the end of this section, we will explain how divided congruences 
might simplify such computations. 
        
Recall \cite[section 4.1]{L2} that if $X$ is a $(U,fr)^2$-manifold then there 
is a decomposition of its
stable tangent bundle $TX=TX^{(0)}\oplus TX^{(1)}$ and we have Chern classes
$c_i^{(j)}\in H^{2i}(X,\Z)$ accordingly ($i\geq 0, j=0,1$).

\begin{theorem}\label{dastheorem}
a) Let $X$ be a $(U,fr)^2$-manifold of dimension $4$ and put $q:=<c_1^{(0)}c_1^{(1)},[X]>\in\Z$. Then $q$ is odd if and only if
the corner of $X$ has Kervaire invariant $1$. If $q$ is even,
then the corner of $X$ is the boundary of a framed manifold.\\
b) Let $X$ be a $(U,fr)^2$-manifold of dimension $8$ and put 
$q:=<c_1^{(0)}(c_1^{(1)3}+c_1^{(1)}c_2^{(1)}+c_3^{(1)})+(c_2^{(0)}+c_1^{(0)2})(c_2^{(1)}+c_1^{(1)2}),[X]>\in\Z$. Then $q$ is odd 
if and only if the corner of $X$ has Kervaire invariant $1$. If $q$ is even,
then the corner of $X$ is the boundary of a framed manifold.
\end{theorem}

\begin{proof} 
If $\sum_{i\geq 0}l_ix^{2^i}$ is the logarithm of the universal $2$-typical formal group law, then (see \cite [Example 4.2.4]{L2})
\[ exp(x)=x-l_1x^2+2l_1^2x^3-(5l_1^3+l_2)x^4\; \; (mod \; x^5)\]
and thus
\[ Q(x):=\frac{x}{exp(x)}=1+l_1x-l_1^2x^2+(2l_1^3+l_2)x^3\; \; (mod \; x^4).\]
For indeterminates $x_i$ of dimension $2$ we set $\Pi:=\prod_i Q(x_i)$.
Denoting by $c_i$ the $i$-th elementary symmetric function in the $x_i$ one gets
(using the definition of the Hazewinkel generators)
\[ \Pi^{(2)}=\frac{v_1}{2}c_1\]
\[ \Pi^{(4)}=\frac{v_1^2}{4}(3c_2-c_1^2) \mbox{ and }\]
\[ \Pi^{(6)}=\frac{v_1^3}{8}(4c_1^3-13c_1c_2+16c_3)+\frac{v_2}{2}(c_1^3-3c_1c_2+3c_3).\]
To prove part a) one has that $\BP_{\Q}^{\otimes 2,4}/(\BP_{\Q}\otimes \Q+\Q\otimes
\BP_{\Q})$ is a one-dimensional $\Q$-vector space 
generated by $v_1\otimes v_1$. Moreover, one checks that the image of
$\BP_4\BP$ is generated (over $\Z_{(2)}$) by $\frac{v_1\otimes v_1}{2}$ and that
$\frac{v_1\otimes v_1}{4}$ is a representative of $\alpha_1^2$.
To see the latter, observe that $\alpha_1$ is represented by $t_1$, hence
\cite[A.1.2.15]{R} $\alpha_1^2$ is represented in the cobar complex 
by $t_1 \otimes t_1=(1 \otimes t_1)(t_1 \otimes 1)$, and then one applies the
description of $\delta^{-1}$ given in section \ref{ratcoeff}.
Using the notations introduced in \cite{L2}, one computes that 

\[ K_{\BP^{<2>}}(TX)^{(4)} = c_1^{(0)}c_1^{(1)}\frac{v_1\otimes v_1}{4}\mbox{ in } \BP_{\Q}^{\otimes 2,4}/(\BP_{\Q}\otimes \Q+\Q\otimes
\BP_{\Q})^{(4)}, \]

hence the image of the corner of $X$ in $\Ext^{2,4}$ is represented by $\frac{q}{2}\cdot\frac{v_1\otimes v_1}{2}=q\cdot\alpha_1^2$. The final assertion follows because the only non-trivial element of
$\pi_2^s$ has Adams-Novikov filtration precisely $2$.\\
For part b), we know that $\Ext^{2,8}$
is generated by $\alpha_1\overline{\alpha}_3$ and $\beta_{2,2}$. As $\beta_{2,2}$ is a permanent cycle in the ANSS
whereas $\alpha_1\overline{\alpha}_3$ is not we know that the image of $X$
in $\Ext^{2,8}$ is a multiple of $\beta_{2,2}$. 
One computes that in the notation of section \ref{ratcoeff}
$\delta(v_2^2)$ is represented by
$z=t_1^4 + v_1^2 t_1^2$ in $C^1(A/2)$.
Hence $\beta_{2,2}=\delta'\delta(v_2^2)$ is represented in 
$\BP_{\Q}^{\otimes 2,8}/(\BP_{\Q}\otimes \Q+\Q\otimes
\BP_{\Q})^{(8)}$ by $-\frac{1}{8}(v_1 \otimes v_1^3) + \frac{5}{16}
(v_1^2 \otimes v_1^2) - \frac{3}{8}(v_1^3 \otimes v_1)$.
Computing enough of the image of $v_1^3t_1,v_1^2t_1^2,v_1^3t_1$ and $t_1^4$ under
$\BP_8\BP\hookrightarrow \BP_8\BP\otimes\Q\simeq \BP_{\Q}^{\otimes 2,8}\longrightarrow 
\BP_{\Q}^{\otimes 2,8}/(\BP_{\Q}\otimes \Q+\Q\otimes
\BP_{\Q})^{(8)}$, one sees that $\frac{v_1^2\otimes v_1^2}{8}\in \imm(\BP_8\BP)$ and
$\beta_{2,2}$ is represented by $\frac{v_1^2\otimes v_1^2}{16}$.
We provide the following argument for general $n$,
for the proof here we need the case $n=3$. 
Observe that by the computations of the previous
sections, the image of $\Ext^{2,2^n}$ in $\BP_{\Q}^{\otimes 2,2^n}/(\BP_{\Q,2^n}\otimes 
\Q+\Q\otimes \BP_{\Q,2^n} + \BP_{2^n}\BP)$ is given by representatives
consisting of summands of the form $v_1^iv_2^k \otimes v_1^j$
with rational coefficients. Moreover, these elements map under
the isomorphism $\phi$ of section \ref{ratcoeff} to polynomials in $v_1,v_2$ and $t_1$.
In other words, the f-invariant in bidegree $(2,2^n)$
factors through the subgroup generated by elements $v_1^iv_2^k \otimes v_1^j$
modulo elements of the form $\phi^{-1}(v_1^iv_2^jt_1^k)$, $1 \otimes v_1^j$
and $v_1^iv_2^j \otimes 1$. Denote this quotient by $B_{2^n}$.
A computation using the results of the previous section shows 
that $\frac{c}{2^{2^{n-1}}}v_1^{2^{n-2}}\otimes v_1^{2^{n-2}}$
is not zero in $B_{2^n}$ for $c$ an odd integer. More precisely,
no element in the relations defining the quotient $B_{2^n}$
contains such a summand. Among the elements in $\Ext^{2,2^n}$,
the image of $\beta_{2^{n-2},2^{n-2}}$ in $B_{2^n}$ contains such a summand
and the other generators exhibited in Proposition \ref{shimo}
do not. Thus we have a well-defined map $\Ext^{2,2^n} \to \Z/2$ given  
by mapping an element to $1$ if it admits a representative 
in $B_{2^n}$ having a summand $\frac{c}{2^{2^{n-1}}}v_1^{2^{n-2}}\otimes v_1^{2^{n-2}}$
with $c$ odd. This map is a projection to $\beta_{2^{n-2},2^{n-2}}$. 
We would like to consider the summands of
\[ K_{\BP^{<2>}}(TX)^{(8)}\equiv \Pi^{(2)}\otimes\Pi^{(6)}+\Pi^{(4)}\otimes\Pi^{(4)}+
\Pi^{(6)}\otimes\Pi^{(2)}\mbox{ mod } (\BP_{\Q}\otimes \Q+ \Q\otimes \BP_{\Q})^{(8)} \]
as elements in $B_{2^3}$. Once this is achieved, 
we have to consider only the summands involving $v_1^2\otimes v_1^2$.
The only summand which is not already given by a representative in 
$B_{2^3}$ is $c_1^{(0)}(c_1^{(1)3}-3c_1^{(1)}c_2^{(1)}+3c_3^{(1)})
\frac{v_1\otimes v_2}{4}$ in $\Pi^{(2)}\otimes\Pi^{(6)}$. One computes 
(for $p=2$ and the Hazewinkel
generators as before) that $\phi^{-1}(t_2)=\frac{1 \otimes v_2}{2}
- \frac{v_2 \otimes 1}{1} + \frac{1 \otimes v_1^3}{4} - \frac{v_1 \otimes 
v_1^2}{8} + \frac{v_1^2 \otimes v_1}{4} -3 \frac{v_1^3 \otimes 1}{8}$.
Hence we have $\phi^{-1}(t_1t_2)=- \frac{v_1 \otimes v_2}{4} - 
\frac{v_1^2 \otimes v_1^2}{16} + ...$, so we have to look at 
the coefficients of $\frac{v_1 \otimes v_2}{4}$ and 
$\frac{v_1^2 \otimes v_1^2}{16}$
which by the above equal $c_1^{(0)}(c_1^{(1)3}-3c_1^{(1)}c_2^{(1)}+3c_3^{(1)})$ 
and $(3c_2^{(0)}-c_1^{(0)2})(3c_2^{(1)}-c_1^{(1)2})$.
We also have that $c_1^{(0)}(c_1^{(1)3}-3c_1^{(1)}c_2^{(1)}+3c_3^{(1)})$ 
equals $c_1^{(0)}(c_1^{(1)3}+c_1^{(1)}c_2^{(1)}+c_3^{(1)})$ and
$(3c_2^{(0)}-c_1^{(0)2})(3c_2^{(1)}-c_1^{(1)2})$
equals $(c_2^{(0)}+c_1^{(0)2})(c_2^{(1)}+c_1^{(1)2})$
modulo $2$ when evaluated on $[X]$. Now the assertion follows 
as in part a).
\end{proof}

Of course, it is possible to do similar but more complicated computations 
for $2^n$-dimensional $(U,fr)^2$-manifolds in case $n \geq 4$. 
We always have a projection to $\Z/2$ looking at the 
power of 2 in the denominator of the coefficient of the summand 
$v_1^{2^{n-2}}\otimes v_1^{2^{n-2}}$. The computation then reduces 
to compute those $\Pi^{(2i)}$ which contribute to this summand.
The diligent reader may thus prove statements of the following form:
The element $h_{n-1}^2$ survives (equivalently: there is a framed
manifold in dimension $2^n-2$ having Kervaire invariant $1$)
if and only if $<F_n(c_i^{(0)},c_j^{(1)}),[X]>$
is odd for a certain explicit polynomial $F_n$.
The main problem in the computation of $F_n$ is to find
representatives in $B_{2^n}$ (that is sums of 
$v_1^iv_2^k \otimes v_1^j$) for elements arising
in the $\Pi^{(2^n-i)}\otimes\Pi^{(i)}$. For $n=3$ this was done using 
$t_1t_2$. In the case $n>3$, it will be necessary to find 
suitable elements in $\BP_{2^n}\BP$ which will involve $t_i$ for
larger $i$ and the computation of $\phi^{-1}$ of these elements.

\medskip

The rest of this section is devoted to the computation of the f-invariant in dimension $2^n$ at the prime $2$.
More generally, we compute the f-invariant of
$\beta_{s2^n,2^n}\in\Ext^{2,(3s-1)2^{n+1}}$ for all $n\ge 0$ and $s\ge 1$ odd.\\
We use the notations of section \ref{div} and those of section \ref{forms}
(in particular the paragraph before Remark \ref{Hecke})
for $p=2$ and $N=3$. We drop the superscripts $\Gamma^{can}(3)$ and $can$ from now 
on if no confusion may arise and write $T:=\frac{a_1-1}{2}\in\D/2\D$ which is an Artin-Schreier generator for the extension $k[a_3]\subseteq V_{1,0}=V_{1,1}\subseteq V_{1,2}$.

\begin{theorem}\label{fbeta}
The image of the f-invariant in $V_{1,2}$ is given by \\
\[ f(\alpha_1\alpha_{2^{n+1}-1})=T\mbox{ for }n\ge 0,\]
\[ f(\beta_s)=1+a_3^s\mbox{ for }s \ge 3 \mbox{ odd }, \]
\[ f(\beta_{s2,2})=1+a_3^{2s}\mbox{ for } s\ge 1\mbox{ odd and}\]
\[ f(\beta_{s2^n,2^n})=(a_3^{2^n}+a_3^{3\cdot 2^{n-2}})^s\mbox{ for }s\ge 1\mbox{ odd and }n\ge 2.\]
\end{theorem}

\begin{proof} For the first line, 
recall that mod $2$ we have $\alpha_t:=\alpha_{t,1}:=\delta(v_1^t)$.
One computes that in the cobar complex $\alpha_1 \alpha_t$
is represented by $t_1 \otimes \frac{1}{2}[(2t_1 + v_1)^t -v_1^t]$, 
use the description of the product in the cobar complex of
\cite[A.1.2.15]{R}. Using the description of $\delta^{-1}$, see section \ref{ratcoeff},
and that $\phi^{-1}$ is a ring isomorphism,
one further computes that $\alpha_1 \alpha_t$ is represented by 
$\frac{1}{4}v_1 \otimes v_1^t$ in the usual quotient of $\BP_{\Q}^{\otimes 2}$, 
c.f. (\ref{inclusion}). Applying the orientation $\alpha$ and $id \otimes q^0$
yields $\frac{1}{4}a_1$ which using (\ref{quot}) and (\ref{mappy})
yields the claim. \\
The second line is a special case of Theorem \ref{betat}
(use Proposition \ref{relation})
and implies the third line because $x_1\equiv x_0^2=v_2^2$ mod $v_1^2$, recall the invariant 
sequences $(2,v_1^{2^n},x_n)$ from section \ref{beta}.\\
The only case requiring a longer computation is $f(\beta_{4,4})$:\\
In the notation of
section \ref{beta} we have $r=1,s=t=4$ and 
\[ x_2=v_2^4-v_1^3v_2^3\in
\HH^{0,24}(\BP/(2,v_1^{4})).\]
This value of $x_2$ follows from the definition in \cite[p. 476]{MRW} 
or \cite[Theorem 5.2.13]{R}) after cancelling all possible multiples of $v_1^4$. One
computes that, in the notation of section \ref{ratcoeff}, $z\in(\overline{\Gamma}/2)^{16}$ is given as
\[ z=t_1^{8}+v_1^{4}t_1^{4}+v_2^{2}t_1(t_1+v_1)+
 v_1v_2t_1^{2}(t_1^{2}+v_1^{2})+v_1^{2}t_1^{3}(t_1+v_1)^3 \]
$$= v_1v_2^{2}t_1
+ v_2^{2}t_1^{2} + v_1^{3}v_2t_1^{2}
+ v_1^{5} t_1^{3} + v_1 v_2
t_1^{4} + v_1^{3} t_1^{5} +
v_1^{2} t_1^{6} + t_1^{8}.$$ 
One then computes that $f(\beta_{4,4})=\nu(w)$, where
$w\in\overline{\Gamma}$ is a lift of $z$ as in section \ref{ratcoeff}, is as claimed, using 
the relation in Proposition \ref{relation} and that
$q^0(a_1)=q^0(a_3)=1$ mod 2 (see Proposition \ref{qexp}).\\
Now the value for $f(\beta_{s4,4})$ for any $s$ follows immediately from the derivation property
of the connecting homomorphism $\delta$, namely $\delta(x^n)=\delta(x)(\sum_{i=0}^{n-1}\eta_R(x)^i\eta_L(x)^{n-1-i})$. 
Alternatively, $f(\beta_{s4,4})$ may be computed directly
for any odd $s$ using that $(q^0 \otimes \id)\eta_L(x_2)=0$ mod 2, $\eta_R$
and $\phi^{-1}$ are ring homomorphisms
and Proposition \ref{relation}.
We obtain $f(\beta_{s2^n,2^n})=f(\beta_{s4,4})^{n-2}$ for all $s$ and $n\ge 3$ because $x_n=x_{n-1}^2$ for all $n\ge 3$.
\end{proof}

Fix some $n\ge 3$. To explain the relevance of the above computation for the problem
of projecting the f-invariant to $\beta_{2^{n-2},2^{n-2}}$ we contemplate the following
diagram.

\begin{equation}\label{diagram}
 \xymatrix{ \Ext^{2,2^n}[\BP] \ar@{^{(}->}[r]^f  \ar@{^{(}->}[rddd]_{\iota} \ar@(ur,ur)[rrr]^{\pi'}
& \Du_{2^{n-1}}\otimes\Q/\Z[2] & & \Z/2 \ar@{^{(}->}[dd] \\
 & \D_{2^{n-1}}\otimes\Q/\Z[2] \ar@{>>}[u]^{(\ref{quot})} \ar@{^{(}->}[r]^(.65){(\ref{mappy})} & \D/2\D & \\
 & \tilde{V_{1,2}} \ar@{^{(}->}[u] \ar@{^{(}->}[r] & V_{1,2} \ar@{^{(}->}[u]_{(\mbox{\small section } \ref{div})} \ar[r]^(.4){\tilde{\pi}} & k=\F_4 \\
 & \tilde{V_{1,2}} \cap k[a_3,T]/(T^2 +T + a_3 + 1)  \ar@{^{(}->}[u] \ar@{^{(}->}[r] &   k[a_3,T]/(T^2 +T + a_3 + 1) \ar@{^{(}->}[u] \ar[ur]^{\pi}  
.} 
\end{equation}

By the results in section \ref{div}, $V_{1,2}$ contains a $k$-free sub-vector
space on the set $\{a_3^iT^j|i\ge 0,j=0,1\}$,
$\pi$ is defined to be the projection to the coefficient of $a_3^{2^{n-2}}$
and $\tilde{\pi}$ is its extension by zero. The map
$\pi'$ is defined to be the projection to the generator $\beta_{2^{n-2},2^{n-2}}$, c.f. Proposition
\ref{shimo}. Theorem \ref{fbeta} determines representatives in $\tilde{V_{1,2}}:=
\D_{2^{n-1}}\otimes\Q/\Z[2]\cap V_{1,2}$ for all generators of $\Ext^{2,2^n}[\BP]$ and thus
defines the map $\iota$. We know that (\ref{diagram}) commutes when $\pi$,
$\tilde{\pi}$ and $\pi'$ are omitted.

\begin{theorem}\label{commutes}
The diagram (\ref{diagram}) is commutative.
\end{theorem}

\begin{proof} By inspection of Proposition \ref{shimo} and Theorem \ref{fbeta}
, the only generator of $\Ext^{2,2^n}[\BP]$ whose f-invariant contains $a_3^{2^{n-2}}$
is the Kervaire element $\beta_{2^{n-2},2^{n-2}}$.
\end{proof}

\begin{cor}\label{detectkervaire}
Let $n\ge 3$ and $X$ a $(U,fr)^2$-manifold of dimension $2^n$.
Then the corner of $X$ has Kervaire invariant one if and only if the f-invariant of
$X$ admits a representative in $\tilde{V_{1,2}}$ which contains the summand $a_3^{2^{n-2}}$.
\end{cor}


\begin{thebibliography}{9999}
\bibitem[APS]{APS} M. Atiyah, V. Patodi, I. Singer, Spectral asymmetry and Riemannian geometry II, Math. Proc. Cambridge Philos. Soc. {\bf 78}  (1975), no. 3, 405--432.
\bibitem[BJM]{BJM} M. Barratt, J. Jones, M. Mahowald,
Relations amongst Toda brackets and the Kervaire invariant in dimension $62$,  
J. London Math. Soc. (2)  {\bf 30}  (1984),  no. 3, 533--550. 
\bibitem[BJM2]{bjm} M. Barratt, J. Jones, M. Mahowald, The Kervaire invariant and the Hopf invariant, in: Algebraic topology (Seattle, Wash., 1985),  135--173, Lecture Notes in Math., {\bf 1286}, Springer, Berlin, 1987.
\bibitem[B]{Br} W. Browder, The Kervaire invariant of framed manifolds
and its generalization, Ann. Math. {\bf 90} (1969), 157-186.
\bibitem[Br]{Bru} R. Bruner, The homotopy theory of $\HH_{\infty}$ ring spectra, in: $\HH_{\infty}$ Ring Spectra and their Applications, Lecture Notes in Math., {\bf 1176}, Springer, Berlin, 1986.
\bibitem[He]{He} E. Hecke, Theorie der Eisensteinschen Reihe h\"oherer
Stufe und ihre Anwendung auf Funktionentheorie und Arithmetik,
Abh. Math. Seminar der Hamb. Univ. {\bf 5} (1927), 199-224
oder Kapitel 24 aus E. Hecke, Mathematische Werke, Vandenhoeck \& Ruprecht,
Göttingen 1959.
\bibitem[HBJ]{H} F. Hirzebruch, T. Berger, R. Jung, Manifolds and modular forms, Aspects of Mathematics, E20, Friedr. Vieweg \& Sohn, Braunschweig, 1992.
\bibitem[HS]{HS} M. Hovey, N. Strickland, Comodules and Landweber exact homology 
theories, Adv. Math. {\bf 192} (2005), 427--456.
\bibitem[K1]{K} N. Katz, Higher congruences between modular forms,
Ann. Math. {\bf 101} (1975), 332--367. 
\bibitem[K2]{Kp} N. Katz, The Eisenstein measure and $p$-adic interpolation,  Amer. J. Math. {\bf 99}  (1977), no. 2, 238--311. 
\bibitem[K3]{K2} N. Katz, $p$-adic properties of modular schemes and modular forms, in: Modular functions of one variable III (Proc. Internat. Summer School, Univ. Antwerp, Antwerp, 1972),  pp. 69--190. Lecture Notes in Mathematics, {\bf 350}, Springer, Berlin, 1973.
\bibitem[KaM]{KaM} N. Katz, B. Mazur, 
Arithmetic moduli of elliptic curves, Annals of Mathematics Studies, {\bf 108},
Princeton University Press, Princeton, NJ, 1985.
\bibitem[KM]{KM} S. Kochman, M. Mahowald,
On the computation of stable stems, Comtemp. Math. {\bf 181} (1995), 299--316.
\bibitem[L1]{L} G. Laures, The topological $q$-expansion principle,
Topology {\bf 38} (1999), 387--425.
\bibitem[L2]{L2} G. Laures, On cobordism of manifolds with corners,
Trans. AMS {\bf 352} (2000), 5667-5688.
\bibitem[MRW]{MRW} H. Miller, D. Ravenel and W. Wilson, 
Periodic phenomena in the Adams-Novikov spectral sequence,
Ann. Math. {\bf 106} (1977), 469--516.
\bibitem[Na]{N} N. Naumann, Comodule categories and the geometry of the stack of formal groups, math.AT/0503308.
\bibitem[Ne]{Ne} J. Neukirch, Algebraic number theory, Grundlehren der Mathematischen Wissenschaften, {\bf 322}, Springer-Verlag, Berlin, 1999.
\bibitem[R]{R} D. Ravenel, Complex cobordism and stable homotopy groups of spheres, Pure and Applied Mathematics, {\bf 121}, Academic Press, Inc., Orlando, FL, 1986. 
\bibitem[Sh]{Sh} K. Shimomura, Novikov's $\Ext^2$ at the prime $2$, Hiroshima Math. 
J. {\bf 11} (1981), 499--513.
\bibitem[Si1]{Si1} J. Silverman, The arithmetic of elliptic curves, Graduate Texts in Mathematics, {\bf 106}, Springer-Verlag, New York, 1986. 
\bibitem[Si2]{Si2} J. Silverman, Advanced topics in the arithmetic of elliptic curves, Graduate Texts in Mathematics, {\bf 151}, Springer-Verlag, New York, 1994. 
\bibitem[St]{St} N. Strickland, Notes on level three structures on elliptic curves,
preprint (2000), available at http://www.shef.ac.uk/personal/n/nps/papers/.
\bibitem[Sw]{Sw} R. Switzer, Algebraic topology---homotopy and homology, Reprint of the 1975 original, Classics in Mathematics, Springer-Verlag, Berlin, 2002.
\bibitem[Wa]{Wa} L. Washington, Introduction to cyclotomic fields, Second edition, Graduate Texts in Mathematics, {\bf 83}, Springer-Verlag, New York, 1997.
\end{thebibliography}
\end{document}